\documentclass{amsart} 
\usepackage{amssymb}
\usepackage{amsmath}
\usepackage{amsthm}
\usepackage{enumerate}
\usepackage{mathrsfs}  
\usepackage{xspace}    
\usepackage{hyperref} 

\newcommand{\abs}[1]{\lvert#1\rvert}
\newcommand{\Adjhm}[1]{\mathbb{B}(#1)}
\newcommand{\Adjlco}[1]{\mathscr{A}_{l,s}^{C_0}(#1)}
\newcommand{\Adjloru}[1]{\mathscr{A}_l^{C_0}(#1)}
\newcommand{\Adjlor}[1]{\mathscr{A}_l(#1)}
\newcommand{\cmcc}{completely contractive\xspace}%
\newcommand{\cmie}{i.\,e.\xspace}
\newcommand{\cmlin}{\operatorname{lin}}%
\newcommand{\cmpmatrix}[1]{\begin{pmatrix}#1\end{pmatrix}}%
\newcommand{\cmsmallpmatrix}[1]{\begin{smallpmatrix}#1\end{smallpmatrix}}%
\newcommand{\cmst}{s.\,t.\xspace}
\newcommand{\cmtimes}{\times}%
\newcommand{\csalgebra}{$C^*$-algebra\xspace}%
\newcommand{\cstern}{$C^*$-}%
\newcommand{\defemph}[1]{\emph{#1\/}}
\newcommand{\dremark}[1]{{[\footnotesize{#1}]}}%
\renewcommand{\dremark}[1]{}%
\newcommand{\dciteda}[1]{\dremark{DA: #1}}
\newcommand{\erzl}[1]{[#1]}
\newcommand{\glqq}{`}%
\newcommand{\grqq}{'}%
\newcommand{\Id}{\operatorname{id}}
\newcommand{\indtHG}{{t \geq 0}}
\newcommand{\indtGr}{{t \in \mbbR}}
\newcommand{\matnull}[1]{#1 \mathrel{\oplus} 0}%
\newcommand{\matoplus}{\oplus}%
\newcommand{\mbbC}{\mathbb{C}}%
\newcommand{\mbbN}{\mathbb{N}}%
\newcommand{\mbbR}{\mathbb{R}}%
\newcommand{\mcalI}{\mathcal{I}}%
\newcommand{\mcalR}{\mathcal{R}}%
\newcommand{\mcalS}{\mathcal{S}}%
\newcommand{\mfrakA}{\mathfrak{A}}%
\newcommand{\mIMl}{\mcalI{}M_\ell}
\newcommand{\mIMls}{\mcalI{}M_\ell^*}
\newcommand{\mLinStet}{L}
\newcommand{\mre}{\mathrm{e}}
\newcommand{\mri}{\mathrm{i}}
\newcommand{\Multcs}[1]{\mathcal{M}(#1)}
\newcommand{\multis}{\cdot}
\newcommand{\Multlco}[1]{\mathscr{M}_l^{C_0}(#1)}%
\newcommand{\Multlor}[1]{{\mathscr{M}_l(#1)}}
\newcommand{\Multwor}[1]{\mcalR(#1)}
\newcommand{\norm}[1]{\ensuremath{\lVert#1\rVert}}
\newcommand{\normbig}[1]{\bigl\lVert#1\bigr\rVert}
\newcommand{\normcb}[1]{\norm{#1}_\text{cb}}
\newcommand{\normcbbig}[1]{\normbig{#1}_\text{cb}}
\newcommand{\normlr}[1]{\left\lVert#1\right\rVert}
\newcommand{\normMl}[2]{\norm{#1}_{\Multlor{#2}}}
\newcommand{\restring}{\ensuremath{\rvert}}   
\newcommand{\setfdg}{\,;\;}
\newcommand{\skalpr}[2]{\langle #1,#2\rangle}
\newcommand{\skalpri}[3]{{\langle #1,#2\rangle}_{#3}}
\newcommand{\sterns}{$\,\!^*$-}
\newcommand{\ternh}[1]{\mathcal{T}(#1)}

\theoremstyle{plain}
\newtheorem{theorem}{Theorem}[section]

\newtheorem{lemma}[theorem]{Lemma}
\newtheorem{proposition}[theorem]{Proposition}

\theoremstyle{definition}
\newtheorem{definitn}[theorem]{Definition}

\newenvironment{enumaequiv}
{\begin{enumerate}[(a)]}%
{\end{enumerate}}
\newenvironment{enumaufz}
{\begin{enumerate}[(i)]}%
{\end{enumerate}}
\newenvironment{smallpmatrix}%
{\left( \begin{smallmatrix}}%
{\end{smallmatrix} \right)}

\begin{document}

\title{Unbounded Multipliers on Operator Spaces}%

%
%

\thanks{The first author was partially supported by the SFB 478 --
Geometrische Strukturen in der Mathematik at the Westf{\"a}lische Wilhelms-Universit{\"a}t M{\"u}nster,
supported by the Deutsche Forschungsgemeinschaft.}

		\author{Hendrik Schlieter}
		\address{Westf\"alische Wilhelms-Universit\"at M\"unster,
               	Mathematisches Institut\\
               	Einsteinstr.~62, D-48149 M\"unster, Germany}
        \email[Hendrik Schlieter]{hschlieter@uni-muenster.de}

		\author{Wend Werner}
		\address{Westf\"alische Wilhelms-Universit\"at M\"unster,
               	Mathematisches Institut\\
               	Einsteinstr.~62, D-48149 M\"unster, Germany}
        \email[Wend Werner]{wwerner@math.uni-muenster.de}

\begin{abstract}
We introduce unbounded multipliers on operator spaces.
These multipliers generalize both,
regular operators on Hilbert \cstern{}modules and
\linebreak 
(bounded) multipliers on operator spaces.
\end{abstract}

\maketitle

\section{Introduction}

Unbounded multipliers play an important role in mathematics,
for example in the spectral theorem for unbounded self-adjoint operators.
In the context of Hilbert \cstern{}modules,
unbounded multipliers are called regular operators.
The definition of such an operator goes back to A.\ Connes, \cite{ConnesThomIsom},
who defines the notion of a self-adjoint regular operator on a \cstern{}algebra.
S.\ L.\ Woronowicz defines in \cite{Woronowicz91} general regular operators
on \cstern{}algebras.
Regular operators are important in the theory of non-compact
quantum groups (\cite{Woronowicz91}, \cite{WoronowiczNapi92OpThInCsAlg})
and can be used to define unbounded Kasparov modules
in Kasparov's bivariant $KK$-theory (\cite{BaajJulg83TheorieBivariante}).

In this article, we define unbounded multipliers on operator spaces.
On a Hilbert \cstern{}module, equipped with a canonical operator space
structure, these multipliers coincide with the regular operators.
Moreover, every bounded left adjointable multiplier on an operator space
belongs to the class of unbounded multipliers.

We will give several characterizations of unbounded multipliers on operator spaces.
For example, it turns out that, roughly, every unbounded multiplier
on an operator space is the restriction of a regular operator
and that an unbounded multiplier can be characterized using the theory of
strongly continuous groups on a Hilbert space.
As an application we will generalize a pertubation result
for regular operators on Hilbert \cstern{}modules to operator spaces.

Most of the results of this paper are from \cite{SchlieterThesis}.

\section{Preliminaries}

For the basic theory of $C_0$-semigroups
we refer the reader to \cite{EngelNagelSemigroups}.
Recall that the generator of a $C_0$-semigroup is
unique, densely defined and closed.

\subsection{Regular operators}

Let $E$, $F$ be Hilbert \cstern{}modules over a \csalgebra{} $\mfrakA$.

\begin{definitn}
Let $A \colon D(A) \subseteq E \to F$ be a densely defined, $\mfrakA$-linear map
where $D(A)$ denotes the domain of $A$.
Then
\[ D(A^*) := \bigl\{ y \in F \setfdg
   \exists z_y \in E\,\, \forall x \in D(A): \skalpr{Ax}{y} = \skalpr{x}{z_y} \bigr\} \]
is a submodule of $F$.
For all $y \in D(A^*)$ the element $z_y$ is uniquely defined and denoted by $A^* y$.
We get an $\mfrakA$-linear, closed map $A^* \colon D(A^*) \subseteq F \to E$,
called the \defemph{adjoint} of $A$, satisfying
$\skalpr{x}{A^*y} = \skalpr{Ax}{y}$ for all $x \in D(A)$ and $y \in D(A^*)$.
\end{definitn}

\begin{definitn}
A \defemph{regular operator} is a densely defined, closed, $\mfrakA$-linear
map $A \colon D(A) \subseteq E \to F$ such that $A^*$ is densely defined and
$1+A^*A$ has dense image.
The set of regular operators is denoted by $\Multwor{E,F}$.
Define $\Multwor{E} := \Multwor{E,E}$.
\end{definitn}

The set of regular operators on a Hilbert space $H$ is equal to
the set of densely defined, closed operators on $H$.
A regular operator $A \in \Multwor{E}$ is called \defemph{self-adjoint} (resp. \defemph{skew-adjoint})
if $A^* = A$ (resp. $A^* = -A)$.
%
%
%
%
For every regular operator $A \in \Multwor{E,F}$ we can define
the adjointable operator $z_A := A\bigl(\left(1+A^*A\right)^{-1}\bigr)^{1/2}$
with $\norm{z_A} \leq 1$.
We denote by $C(\mbbR)$ the set of continuous functions from $\mbbR$ to $\mbbC$.

For self-adjoint regular operators we have the following functional calculus:

\begin{theorem}[{\cite[Theorem 10.9]{LanceHmod}}]
Let $A \in \Multwor{E}$ be self-adjoint.
Let $\iota$ be the canonical inclusion of $\mbbR$ into $\mbbC$ and
$f \colon \mbbR \to \mbbC, \lambda \mapsto \lambda \left(1+\lambda^2\right)^{-1/2}$.
Then there exists a unique \sterns{}homomorphism $\varphi_A \colon C(\mbbR) \to \Multwor{E}$ such that
$\varphi_A(\iota) = A$ and $\varphi_A(f) = z_A$.
\dremark{Hier ist u. U. nicht klar, was genau mit \sterns{}Hom. gemeint ist,
  denn $\Multwor{E}$ ist keine \sterns{}Algebra.}%
\end{theorem}

\subsection{Stone's theorem for Hilbert \cstern{}modules}

Let $E$ be a Hilbert \cstern{}module over a \csalgebra $\mfrakA$.
Denote by $\Adjhm{E}$ the set of adjointable operators on $E$.
The following two results are generalizations of
\cite[Theorem 2.1 and Proposition~2.2]{HollevoetEtAl92StonesTh}
from \cstern{}algebras to Hilbert \cstern{}modules.
With some adjustments, the proof of the first statement
can be transferred to the Hilbert \cstern{}module case.
For details see \cite{SchlieterThesis}.

\begin{theorem}\label{erzeugerUnitaereGrWor}
Let $A \in \Multwor{E}$ be self-adjoint.
Let $e_t \colon \mbbR \to \mbbC, \lambda \mapsto \exp(\mri \lambda t)$, and
$U_t := \exp(\mri tA) := \varphi_A(e_t)$ for all $\indtGr$.
Then $(U_t)_\indtGr$ is a $C_0$-group on $E$ with generator $\mri A$
such that $U_t$ is a unitary element of $\Adjhm{E}$ for all $\indtGr$.
\dciteda{1.27, 1.30}
\end{theorem}



\begin{theorem}[Stone's theorem]\label{HollevoetTh21}
Let $(U_t)_\indtGr$ be a $C_0$-group on $E$
such that $U_t$ is a unitary element of $\Adjhm{E}$ for all $\indtGr$.
\begin{enumaufz}
\item There exists a self-adjoint  $A \in \Multwor{E}$ such that
$U_t = \exp(\mri tA)$  for all $\indtGr$.
\item If $(U_t)_\indtGr$ is norm continuous, then $A \in \Adjhm{E}$.\dciteda{1.29}
\end{enumaufz}
\end{theorem}

\begin{proof}
We will only prove (i).
Denote by $C^*(\mbbR)$ the group \cstern{}algebra of the locally compact group $\mbbR$
and by $\hat{\mbbR}$ the dual group of $\mbbR$.
We have $\hat{\mbbR} \cong \mbbR$.
Let $C_0(\mbbR)$ (resp. $C_c(\mbbR)$) denote the set of
continuous functions on $\mbbR$ vanishing at infinity (resp. with compact support).
For all $f \in C_c(\mbbR)$
\[ \alpha(f) = \int_\mbbR f(t) U_t\, dt \]
defines an operator in $\Adjhm{E}$
and the map 
$\alpha$ can be extended to \sterns{}homomorphism
$\check{\alpha} \colon \Multcs{C^*(\mbbR)} \to \Adjhm{E}$ (\cite[Proposition C.17]{RaeburnWilliams98MoritaEq})
where $\Multcs{C^*(\mbbR)}$ denotes
the algebra of double centralizers of the algebra $C^*(\mbbR)$.

Let $t \in \mbbR$.
Define
\[ \lambda_t \colon L^1(\mbbR) \to L^1(\mbbR), f \mapsto (s \mapsto f(s-t) ). \]
Using the translation invariance of the Lebesgue measure,
we get $\lambda_t \in \Multcs{L^1(\mbbR)}$.
By approximating any $f \in C^*(\mbbR)$ by elements of $L^1(\mbbR)$,
we can define the double centralizer 
$\hat{\lambda}_t \in \Multcs{C^*(\mbbR)}$ using $\lambda_t$.
For all $f \in C_c(\mbbR)$ and $x \in E$ we obtain
\begin{align*}
   \check{\alpha}(\hat{\lambda}_t) \alpha(f) x
=  \alpha(\lambda_t f) x
=  U_t \alpha(f) x
\end{align*}
and thus
\begin{equation}\label{eqalphaltEqUt2}
\check{\alpha}(\hat{\lambda}_t) = U_t.
\end{equation}

The Fourier transform
$\mathscr{F} \colon L^1(\mbbR) \to C_0(\hat{\mbbR}) \cong C_0(\mbbR)$
can be extended to a \sterns{}iso\-mor\-phism
$\hat{\mathscr{F}} \colon C^*(\mbbR) \to C_0(\mbbR)$
(\cite[Example C.20]{RaeburnWilliams98MoritaEq}).
Denote by $C_b(\mbbR)$ the set of bounded continuous functions on $\mbbR$.
Because there is an embedding of $C_0(\mbbR)$ into $\Multcs{C_0(\mbbR)} \cong C_b(\mbbR)$,
we can regard $\hat{\mathscr{F}}$ as a map from $C^*(\mbbR)$ to $\Multcs{C_0(\mbbR)}$.
Using \cite[Proposition 2.5]{LanceHmod},
we can extend this nondegenerate \sterns{}homomorphism to a \sterns{}homomorphism
$\check{\mathscr{F}} \colon \Multcs{C^*(\mbbR)} \to \Multcs{C_0(\mbbR)}$
and conclude
\begin{equation}\label{eqFlambdatEqet}
\check{\mathscr{F}}(\hat{\lambda}_t) = e_t.
\end{equation}

The nondegenerate \sterns{}homomorphism
$\beta := \hat{\alpha} \circ \hat{\mathscr{F}}^{-1} \colon C_0(\mbbR) \to \Adjhm{E}$
extends to a \sterns{}homomorphism
$\check{\beta} \colon \Multcs{C_0(\mbbR)} \to \Adjhm{E}$.
Recall that $\iota$ is the canonical inclusion of $\mbbR$ into $\mbbC$.
Using $\iota \in C(\mbbR) \cong \Multwor{C_0(\mbbR)}$ (\cite[Example 2]{Woronowicz91}),
it follows that $A := \beta(\iota)$ is a self-adjoint regular operator on $E$.
One can show that $\varphi_A\restring_{C_b(\mbbR)} = \check{\beta}$
and $\check{\beta} \circ \check{\mathscr{F}} = \check{\alpha}$.
Together with \eqref{eqFlambdatEqet} und \eqref{eqalphaltEqUt2}, we obtain
\[ \exp(\mri tA)
=  \varphi_A(e_t)
=  \check{\beta}(e_t)
=  \check{\beta}\left(\check{\mathscr{F}}(\hat{\lambda}_t)\right)
=  \check{\alpha}(\hat{\lambda}_t)
=  U_t.  \qedhere \]
\end{proof}

\subsection{Operator spaces}

Let $H$ be a Hilbert space.
Denote by $\mLinStet(H)$ the algebra of linear, continuous operators on $H$.
A \defemph{(concrete) operator space} is a closed subspace $X$ of $\mLinStet(H)$.
Let $M_{m,n}(X)$ denote the space of $m \times n$ matrices with entries in $X$,
and set $M_n(X) := M_{n,n}(X)$.
We may view the space $M_n(X)$ as a subspace of $\mLinStet(H^n)$.
Thus this matrix space has a natural norm for each $n \in \mbbN$.
Note that operator spaces can be characterized intrinsically.


Let $X$, $Y$ be operator spaces and $\alpha \colon X \to Y$ linear.
Define
\[ \alpha_n := \alpha^{(n)} : M_n(X) \rightarrow M_n(Y), x \mapsto (\alpha(x_{ij}))_{i,j}, \]
called the $n$th \defemph{amplification} of $\alpha$.
The map $\alpha$ is called \defemph{completely contractive}
(resp. \defemph{completely isometric}) if $\alpha_n$ is contractive (resp. isometric)
for all $n \in \mbbN$.
Furthermore, $\alpha$ is called \defemph{completely bounded}
if $\normcb{\alpha} := \sup_{ n \in \mbbN } \norm{\alpha_n} < \infty$.


In the following, let $X$ be an operator space.

The operator space $X$ is called \defemph{unital}
if it has a distinguished element $e_X$
such that there exist a \cstern{}algebra $\mfrakA$ with unit $e_\mfrakA$ and
a completely isometric map $\eta \colon X \to \mfrakA$
satisfying $\eta(e_X) = e_\mfrakA$.
Observe that unital operator spaces can be characterized intrinsically
(\cite{BlecherNeal08MetricChar}, \cite{HuangNg08AbstractCharUnitalOS}).
\dremark{Am Ende auf Ver"offentlichung "uberpr"ufen}%

Recall that $X$ is called \defemph{injective}
if for every operator space $Y \subseteq Z$
and every completely contractive map $\varphi \colon Y \to X$
there is a completely contractive map $\Phi \colon Z \to X$
such that $\Phi\restring_Y = \varphi$.
For every operator space $X$ there exists a unique \defemph{injective envelope}
$(I(X),j_X)$, \cmie, a completely isometric map $j_X \colon X \to I(X)$
and an injective operator space $I(X)$ such that
for every injective subspace $Y$ of $I(X)$
satisfying $j_X(X) \subseteq Y$ we have $Y = I(X)$.

Let $X \subseteq \mLinStet(H)$ be an operator space.
We can embed $X$ into an \defemph{operator system},
\cmie, a self-adjoint closed subspace of $\mLinStet(H)$ containing $\Id_H$.
To accomplish this we form the Paulsen system
\[ \mcalS(X) :=
\begin{pmatrix} \mbbC \Id_H & X \\ X^* & \mbbC \Id_H \end{pmatrix}
\subseteq M_2(\mLinStet(H)). \]
There exists a completely contractive map
$\Phi \colon M_2(\mLinStet(H)) \to M_2(\mLinStet(H))$
such that the image of $\Phi$ is an injective envelope of the operator space $\mcalS(X)$.
Moreover, $I(\mcalS(X))$ is a unital \cstern{}algebra
with multiplication $x \multis y := \Phi(xy)$.
Let $p_1 := \Id_H \mathrel{\oplus} 0$ and $p_2 := 0 \oplus \Id_H$
be the canonical projections in $I(\mcalS(X))$.
Then $p_1$ and $p_2$ are orthogonal and corner preserving.
Therefore, with respect to $p_1$ and $p_2$,
we may decompose $I(\mcalS(X))$ as
\[ I(\mcalS(X))
= \begin{pmatrix} I_{11}(X) &  I_{12}(X) \\ I_{12}(X)^* & I_{22}(X) \end{pmatrix} \]
and have $I_{12}(X) = I(X)$.
\smallskip


The following operator spaces will play a crucial role in this article.

If $\mfrakA$ is a \cstern{}algebra, then the \sterns{}algebra
$M_n(\mfrakA)$ has a unique norm such that $M_n(\mfrakA)$ is a \cstern{}algebra.
With respect to these norms, $\mfrakA$ is an operator space.

Let $X$ be an operator space and $m, n \in \mbbN$.
We can regard $M_{m,n}(X)$ as a subspace of $M_p(X)$ where $p := \max\{m,n\}$.
Thus $M_{m,n}(X)$ is an operator space.
The set $C_n(X) := M_{n,1}(X)$ is called \defemph{column operator space}.

Furthermore, every Hilbert \cstern{}module $E$ over a \cstern{}algebra $\mfrakA$
carries the structure of an operator space.
For every $n \in \mbbN$ and $x \in M_n(E)$
\begin{equation}\label{eqHCsModAsOS}
\norm{x}_n = \normlr{ \left( \sum_{k=1}^n \skalpr{x_{ki}}{x_{kj}} \right)_{i,j} }^{1/2}
\end{equation}
defines a norm such that $E$ becomes an operator space.
\dremark{Frage: Zusammenhang zu Zeilen-/Spaltenoperatorraum?}

\subsection{Multipliers on operator spaces}

Let us recall some facts about multipliers on operator spaces.
These multipliers generalize multipliers on \cstern{}algebras.

\begin{definitn}\label{defLinksmultOR}
Let $T \colon X \to X$ be linear.
\begin{enumaufz}
\item $T$ is called \defemph{(left) multiplier} on $X$
if there exist a Hilbert space $H$ together with
a completely isometric map $\eta \colon X \to \mLinStet(H)$
and $a \in \mLinStet(H)$ such that $\eta(Tx) = a\eta(x)$
for all  $x \in X$.
The set of all multipliers on $X$ is denoted by $\Multlor{X}$.

\item For all $T \in \Multlor{X}$ define the multiplier norm
\begin{align*}
   \normMl{T}{X}
   &:=  \inf\bigl\{ \norm{a} \setfdg \text{there exist a Hilbert space } H,\text{ a completely} \\
   &\phantom{:=  \inf\bigl\{ \norm{a} \setfdg} \,\,\,\text{isometric map } \eta \colon X \to \mLinStet(H) \text{ and } a \in \mLinStet(H) \text{ \cmst }\\
   &\phantom{:=  \inf\bigl\{ \norm{a} \setfdg} \,\,\,\forall x \in X: \eta(Tx) = a\eta(x) \bigr\}.
\end{align*}

\item $T$ is called \defemph{left adjointable multiplier} on $X$
if there exist a Hilbert space $H$, a completely isometric map $\eta \colon X \to \mLinStet(H)$
and a map $S \colon X \to X$ such that
$\eta(Tx)^*\eta(y) = \eta(x)^*\eta(Sy)$ for all $x,y \in X$.
The set of all left adjointable multipliers on $X$ is denoted by $\Adjlor{X}$.
\end{enumaufz}
\end{definitn}

We have $\Multlor{X} \subseteq \Adjlor{X}$.
Moreover, $\Adjlor{X}$ is a unital \cstern{}algebra.
The multiplier algebra of a \cstern{}algebra $\mfrakA$
is isomorphic to $\Adjlor{\mfrakA}$.
If $E$ is a Hilbert \cstern{}module, then
\begin{equation}\label{eqBeEqAlE}
\Adjlor{E} = \Adjhm{E}
\end{equation}
(\cite[Corolarry 8.4.2]{BlecherLeMerdy04OpAlg}).
Define
\begin{align*}
\mIMl(X) &:= \{ a \in I_{11}(X) \setfdg a \multis X \subseteq X \}
\quad\text{and}  \\
\mIMls(X) &:= \mIMl(X) \cap \mIMl(X)^*.
\end{align*}
One can show that $\mIMl(X)$ (resp. $\mIMls(X)$) is isometrically isomorphic to
$\Multlor{X}$ (resp. $\Adjlor{X}$).

Left multipliers are characterized intrinsically in the following result:

\begin{theorem}[{\cite[Corollary 3.13]{WWerner03}, \cite[Theorem 1.1]{BlecherEffrosZarikian02}}]\label{charMultORmitLambda}\label{charMultOR}
Let $T \colon X \to X$ be linear and $\lambda > 0$.
The following assertions are equivalent:
\begin{enumaequiv}
\item $T \in \Multlor{X}$ satisfies $\normMl{T}{X} \leq \frac{1}{\lambda}$.
\item The map
\[ \tau_{\lambda T} \colon C_2(X) \to C_2(X), \binom{x}{y} \mapsto \binom{\lambda Tx}{y}, \]
  is \cmcc.
\dremark{Literaturangabe pr"ufen}%
\end{enumaequiv}
\end{theorem}

\begin{proposition}[{\cite[Proposition 1.7.6. and 1.7.8]{Zarikian01Thesis}}]\label{charAlXZar}
Let $U \colon X \to X$ be linear.
The following assertions are equivalent:
\begin{enumaequiv}
\item $U$ is a unitary element of $\Adjlor{X}$.
\item The map
\[ \tau_U \colon C_2(X) \to C_2(X), \binom{x}{y} \mapsto \binom{Ux}{y}, \]
  is a completely isometric isomorphism.
\item There exists a unitary $u \in \mIMls(X)$ such that
$j(Ux) = u \multis j(x)$ for all $x \in X$.
\dremark{Dies steht z. T. auch in: \cite{BlecherEffrosZarikian02}, Cor. 4.8.}
\end{enumaequiv}
\end{proposition}

\section{Unbounded multipliers on operator spaces}

In this chapter, we will introduce three different definitions
of unbounded multipliers on an operator space.

\subsection{Unbounded skew-adjoint multipliers}

An operator space has less structure than a Hilbert \cstern{}module.
Therefore, it is not obvious, how to define an unbounded multiplier
on an operator space.
We know that the set $\Adjhm{E}$ of adjointable operators
is isomorphic to the set $\Adjlor{E}$,
where $E$ is equipped with the canonical operator space structure of equation
\eqref{eqHCsModAsOS}.
The fact that skew-adjoint regular operators can be
characterized using $C_0$-groups
(Theorem \ref{HollevoetTh21} and \ref{erzeugerUnitaereGrWor})
motivates the following definition:

\begin{definitn}
\begin{enumaufz}
\item An operator $A \colon D(A) \subseteq X \to X$ is called \defemph{$C_0$-left multiplier}
if $A$ generates a $C_0$-semigroup $(T_t)_\indtHG$ on $X$ satisfying $T_t \in \Multlor{X}$
for all $\indtHG$.

\item An operator $A \colon D(A) \subseteq X \to X$ is called \defemph{unbounded skew-adjoint multiplier}
if $A$ generates a $C_0$-group $(U_t)_\indtGr$ on $X$ satisfying $U_t \in \Adjlor{X}$
for all $\indtGr$.

\item The set of all $C_0$-left multipliers (resp. unbounded skew-adjoint multipliers)
on $X$ is denoted by $\Multlco{X}$ (resp. $\Adjlco{X}$).
\end{enumaufz}
\end{definitn}

Using the fact that $t \mapsto \mre^{tA}$ ($A \in \mLinStet(X)$)
defines a norm-continuous semigroup with generator $A$, we get:

\begin{proposition}\label{MultlorSubseteqMultlco}\label{AdjlorSubseteqAdjlco}
\begin{enumaufz}
\item $\Multlor{X} = \{ A \in \Multlco{X} \setfdg D(A) = X \}$,
\item $\{ A \in \Adjlor{X} \setfdg A \text{ skew-adjoint} \}
  = \{ A \in \Adjlco{X} \setfdg D(A) = X \}$.
\end{enumaufz}
\end{proposition}

In a unital \cstern{}algebra $\mfrakA$ the set of (bounded) multipliers
is equal to $\Multwor{\mfrakA}$.
A similar result holds for operator spaces:

\begin{proposition}\label{AdjlorsaEqAdjlco}
If $X$ is unital, we have
\begin{enumaufz}
\item $\Multlor{X} = \Multlco{X}$ and
\item $\bigl\{ A \in \Adjlor{X} \setfdg A \text{ skew-adjoint} \bigr\} = \Adjlco{X}$.\dciteda{3.8}
\end{enumaufz}
\end{proposition}

\begin{proof}
We sketch the proof of the inclusion \glqq$\supseteq$\grqq{} of (i).
Suppose $A \in \Multlco{X}$ generates the $C_0$-semigroup $(T_t)_\indtHG$.
There exists a Hilbert space $H$
such that $X \subseteq \mLinStet(H)$ and $e_X = \Id_H$.
Because $X$ is unital, we have $\normcb{S} = \normMl{S}{X}$
for all $S \in \Multlor{X}$ (\cite[p.\ 20]{Blecher01ShilovBoundary}).
One can see that
\begin{align*}
   \norm{T_t - \Id_X}
&\leq  \normcb{T_t - \Id_X}
=  \normMl{T_t - \Id_X}{X}
= \norm{ T_t(e_X) - e_X }
\end{align*}
converges to $0$ for $t \to 0$.
Thus $(T_t)_\indtHG$ is norm-continuous.
Using $A = (t \mapsto T_t')(0)$, the assertion follows.
\dremark{Bew.: Sehr knapp}%
\end{proof}

\begin{theorem}\label{MultworInAdjlco}
For a Hilbert \cstern{}module $E$, we have
\[ \left\{ A \in \Multwor{E} \setfdg A \text{ skew-adjoint} \right\} = \Adjlco{E}. \]
\end{theorem}

\begin{proof}
Let $A \in \Multwor{E}$ be skew-adjoint.
Set $h := -\mri A$.
Then $U_t := \exp(\mri th)$ defines a $C_0$-group on $E$
with $U_t \in \Adjhm{E} = \Adjlor{E}$ unitary
and with generator $\mri h = A$ due to Theorem \ref{erzeugerUnitaereGrWor}
and equation \eqref{eqBeEqAlE}.

The other inclusion follows with Theorem \ref{HollevoetTh21}, Theorem \ref{erzeugerUnitaereGrWor}
and the uniqueness of the generator of a $C_0$-group.
\end{proof}

\subsection{Unbounded multipliers}

In this section we will generalize the definition of
regular operators on Hilbert \cstern{}modules to operator spaces.

Let $A_i \colon D(A_i) \subseteq X \to X$ be an operator for all $i \in \{1, \dots, 4\}$.
Define
\begin{align*}
 \cmpmatrix{A_1 & A_2 \\ A_3 & A_4} \colon
&\,\,(D(A_1) \cap D(A_3)) \cmtimes (D(A_2) \cap D(A_4)) \subseteq C_2(X)  \to  C_2(X),  \\
&\,\,(x,y) \mapsto (A_1 x + A_2 y, A_3 x + A_4 y).
\end{align*}

\begin{definitn}
A map $A \colon D(A) \subseteq X \to X$ is called
\defemph{unbounded multiplier}
if there exists $B \colon D(S) \subseteq X \to X$ satisfying
$\mri \begin{pmatrix} 0 & A \\ B & 0 \end{pmatrix} \in \Adjlco{C_2(X)}$.
The set of all unbounded multipliers on $X$ is denoted by $\Adjloru{X}$.
\end{definitn}

\begin{proposition}\label{AdjlorSubseteqAdjloru}\label{XunitalAdjlorEqAdjloru}
\begin{enumaufz}
\item $\Adjlor{X} \subseteq \Adjloru{X}$.
\item If $X$ is unital, we have $\Adjlor{X} = \Adjloru{X}$.
\end{enumaufz}
\end{proposition}

\begin{proof}
(i) follows using Proposition \ref{AdjlorSubseteqAdjlco}.
Part (ii) is deduced using Proposition \ref{AdjlorsaEqAdjlco} and the following lemma.
\end{proof}

The following is straightforward using \cite[Corollary 1.3]{BlecherPaulsen00MultOfOpSpaces}:

\begin{lemma}\label{NullTS0AdjFolgtTAdj}
Let $A \colon D(A) \subseteq X \to X$ and $B \colon D(B) \subseteq X \to X$ be operators
satisfying $\hat{A} := \cmpmatrix{0 & A \\ B & 0} \in \Adjlor{C_2(X)}$ and $(\hat{A})^* = \hat{A}$.
Then $A, B \in \Adjlor{X}$ and $B^* = A$.\dciteda{3.20}
\end{lemma}

The following important theorem shows
that on a Hilbert \cstern{}module unbounded multipliers
coincide with regular operators:

\begin{theorem}\label{zshgRegOpUnbeschrMult}
Let $E$ be a Hilbert \cstern{}module.
\begin{enumaufz}
\item $\Adjloru{E} = \Multwor{E}$.
\item Let $A \in \Adjloru{E}$, and
suppose there exists $B \colon D(B) \subseteq E \to E$ such that
$\mri \cmpmatrix{0 & A \\ B & 0} \in \Adjlco{C_2(E)}$.
Then $A^* = B \in \Multwor{E}$.\dciteda{3.21}
\end{enumaufz}
\end{theorem}

\begin{proof}
We first show (ii) and the inclusion \glqq$\subseteq$\grqq{} of (i).
Let $A \in \Adjlco{E}$ and $\hat{A} := \cmpmatrix{0 & A \\ B & 0}$.
The direct sum $E \oplus E$ is, equipped with
the inner product
\[ \skalpri{(x_1,y_1)}{(x_2,y_2)}{E \oplus E} := \skalpri{x_1}{x_2}{E} + \skalpri{y_1}{y_2}{E}, \]
a Hilbert \cstern{}module.
One can show
that $E \oplus E$ is isometrically isomorphic
to the column operator space $C_2(E)$ of $E$ viewed as an operator space.
We obtain with Theorem \ref{MultworInAdjlco}:
$\hat{A} \in \Adjloru{C_2(E)} \cong \Multwor{E \oplus E}$ is self-adjoint,
thus
\[ \hat{A} = \bigl(\hat{A}\bigr)^* = \cmpmatrix{0 & B^* \\ A^* & 0}. \]
It follows $B = A^* \in \Multwor{E}$.

The remaining part follows using Theorem \ref{MultworInAdjlco}.
\end{proof}

\subsection{Characterizations of unbounded multipliers}

Our next goal is to characterize unbounded multipliers.
Let us recall the following two definitions:

\begin{definitn}
Let $H$, $K$ be Hilbert spaces.
A \defemph{ternary ring of operators} (\defemph{TRO})
is a closed subspace $Z$ of $\mLinStet(H,K)$
such that $Z Z^* Z \subseteq Z$.
A \defemph{subtriple} of a TRO $Z$ is a closed subspace $Y$ of $Z$
satisfying $Y Y^* Y \subseteq Y$.
\end{definitn}

Every TRO $Z$ is a full Hilbert \cstern{}module over the \cstern{}algebra $Z^* Z$.
Conversely, every Hilbert \cstern{}module can be represented faithfully
as a TRO.
We can view the injective envelope $I(X)$ as a TRO.

\begin{definitn}
The triple envelope $\ternh{X}$ of $X$ is the smallest subtriple
of $I(X)$ containing $X$.
\end{definitn}

We have $\ternh{X} = \overline{\cmlin}
  \{ x_1 \multis x_2^* \multis x_3 \multis x_4^* \cdots x_{2n+1} \setfdg n \geq 0, x_1, \dots , x_{2n+1} \in j(X) \}$.

The following is straightforward from \cite[p. 341]{Harris81GeneralizationOfCsAlg}:

\begin{lemma}\label{TXEmbeddingExists}
There exist a Hilbert space $K$,
a unital injective \sterns{}representation $\pi \colon I(\mcalS(X)) \to \mLinStet(K)$
and closed subspaces $K_1$, $K_2$ of $K$ such that
\begin{enumaufz}
\item $\erzl{\pi(\ternh{X})K_1} := \cmlin\{ x y \setfdg x \in \pi(\ternh{X}), y \in K_1\}$ is dense in $K_2$,
\item $\erzl{\pi(\ternh{X})K_1^\perp} = \{0\}$ and
\item $\erzl{\pi(j(X)) K} \subseteq K_2$.\dciteda{3.62}
\end{enumaufz}
\end{lemma}

\begin{definitn}
The quadrupel $(\pi, K, K_1, K_2)$, introduced in Lemma \ref{TXEmbeddingExists},
is called a \defemph{$\ternh{X}$-embedding} of $X$.
\end{definitn}

In the following theorem it is shown that
an unbounded multiplier can be characterized
intrinsically (part (b))
or with the help of a $C_0$-group on a Hilbert space (part (d))
and that an unbounded multiplier is essentially
the restriction of a regular operator on the Hilbert \cstern{}module $\ternh{X}$ (part (c)).
\pagebreak 

\begin{theorem}\label{thMainChar}
Let $(\pi,K,K_1,K_2)$ be a $\ternh{C_2(X)}$-embedding of $C_2(X)$, and set
$\eta := \pi \circ j_{C_2(X)} \colon C_2(X) \to \mLinStet(K)$.
Let $A \colon D(A) \subseteq X \to X$ be linear.
The following assertions are equivalent:
\begin{enumaequiv}
\item $A \in \Adjloru{X}$.

\item There exists a map $B \colon D(B) \subseteq X \to X$ such that
  $\mri \cmpmatrix{ 0 & A \\ B & 0 } \matoplus \cmpmatrix{ 0 & 0 \\ 0 & 0 }$
  generates a \cmcc $C_0$-group on $C_2(C_2(X)) \cong C_4(X)$.\dciteda{3.29}

\item There exists $B \in \Multwor{\ternh{X}}$ with the following properties:
\begin{enumaufz}
\item $j \circ A = B \circ j\restring_{D(A)}$ and
\item there exist $\lambda, \mu > 0$ such that the restriction of
\[ \lambda - \mri\cmpmatrix{0 & B \\ B^* & 0} \quad\text{resp.}\quad
   \mu + \mri\cmpmatrix{0 & B \\ B^* & 0} \]
  to $j_{C_2(X)}(C_2(X))$ is surjective onto $j_{C_2(X)}(C_2(X))$.\dciteda{3.59}
\end{enumaufz}

\item  There exist a map $B \colon D(B) \subseteq X \to X$
  such that $\tilde{A} := \mri \cmsmallpmatrix{ 0 & A \\ B & 0}$
  is densely defined with $\rho(\tilde{A}) \neq \emptyset$
  and a $C_0$-group $(b_t)_\indtGr$
  consisting of unitary elements of $K_2$ such that:
\begin{enumaufz}
\item $L_{b_t} \colon \eta(C_2(X)) \to \eta(C_2(X)), y \mapsto b_t y,$
  and $L_{b^*_t}$ define $C_0$-groups on the space 
  $\eta(C_2(X))$ and
\item for the generator $C$ of $(L_{b_t})_\indtGr$, we have $\eta(D(\tilde{A})) \subseteq D(C)$.

%
\dciteda{3.67}%
\end{enumaufz}
\end{enumaequiv}
\end{theorem}

Note that part (c).(i) of the above theorem implies $j(D(A)) \subseteq D(B)$.

A first step in the proof of the above characterization is the following result:
\dremark{Ziel: Den obigen Satz beweisen. Zun"achst:}

\begin{theorem}\label{A0ErzvkC0GrFolgt}\label{charMultlcomatnull}\label{charAdjlcomatnull}
Let $A \colon D(A) \subseteq X \to X$.
Then $A \in \Adjlco{X}$ if and only if
$\matnull{A}$ generates a \cmcc $C_0$-group on $C_2(X)$.\dciteda{3.26}

\end{theorem}

\begin{proof}
This follows with Proposition \ref{charAlXZar}.
\end{proof}

%


With this result we have proven the equivalence of (a) and (b) in Theorem \ref{thMainChar}.
Thus unbounded multipliers can be characterized
using generators of completely contractive $C_0$-groups.
Such generators are analyzed in the following two theorems,
which generalize well-known results from the theory
of contractive $C_0$-semigroups to the setting of \cmcc (c.\,c.{}) $C_0$-semigroups
on an operator space~$X$.
With appropriate adjustments, one can prove both theorems.


\begin{theorem}[Theorem of Hille-Yosida for c.\,c.\ $C_0$-se\-mi\-groups]\label{satzHilleYosvkHG}
An operator $A \colon D(A) \subseteq X \to X$ is the generator of a \cmcc $C_0$-semigroup
if and only if $A$ is densely defined and closed,
$\mbbR_{>0}$ is a subset of the resolvent set $\rho(A)$ of $A$ and
\begin{equation*}
\normcb{ \lambda (\lambda - A)^{-1} } \leq 1
   \qquad \text{for all } \lambda > 0.
\end{equation*}
\end{theorem}

An operator $A \colon D(A) \subseteq X \to X$ is called \defemph{completely dissipative}
if $A_n$ is dissipative for all $n \in \mbbN$, \cmie,
\[ \normlr{(\lambda-A_n)x} \geq \lambda \norm{x}
   \quad\text{for all } \lambda > 0 \text{ and } x \in D(A_n). \]

\begin{theorem}[Theorem of Lumer-Phillips for c.\,c.\ $C_0$-semigroups]\label{SatzLumerPhillipsOR}
\rule{1pt}{0pt}\\
Let 
$A \colon D(A) \subseteq X \to X$ be linear and densely defined.
Then $A$ generates a \cmcc $C_0$-semigroup
if and only if $A$ is completely dissipative and
there exists a $\lambda > 0$ such that
$\lambda - A$ is surjective onto $X$.
\end{theorem}

Using the two theorems above together with
the equivalence of (a) and (b) in Theorem \ref{thMainChar},
we obtain:

\begin{theorem}
Let $A \colon D(A) \subseteq X \to X$ be linear.
The following assertions are equivalent:
\begin{enumaequiv}
\item $A \in \Adjloru{X}$.

\item There exists $B \colon D(B) \subseteq X \to X$ such that
  $\tilde{A} := \mri \cmsmallpmatrix{0 & A \\ B & 0}$ is densely defined and closed
  satisfying $\mbbR\setminus\{0\} \subseteq \rho(\tilde{A})$ and
\[ \normcbbig{\abs{\lambda} \, \bigl(\lambda - (\matnull{\tilde{A}})\bigr)^{-1} } \leq 1 \qquad\text{for all } \lambda \in \mbbR\setminus\{0\}. \]

\item There exists $B \colon D(B) \subseteq X \to X$ with the following properties:
\begin{enumaufz}
\item $\tilde{A} := \mri \cmsmallpmatrix{0 & A \\ B & 0}$ is densely defined,
\item $\matnull{\tilde{A}}$ and $-\matnull{\tilde{A}}$ are completely dissipative and
\item there exist $\lambda, \mu >0$ such that
  $\lambda - \tilde{A}$ and $\mu + \tilde{A}$ are surjective onto $C_2(X)$.
\end{enumaufz}
\end{enumaequiv}
\end{theorem}


The following result shows that we can lift unbounded multipliers from $X$ to the
ternary envelope $\ternh{X}$:

\begin{proposition}\label{grHochliftenAufTX}\label{AdjloruHochliftenAufTX}
Let $A \in \Adjlco{X}$ (resp. $A \in \Adjloru{X}$).
There exists a unique $B \in \Adjlco{\ternh{X}}$ (resp. $B \in \Adjloru{\ternh{X}}$) such that
$j \circ A = B \circ j\restring_{D(A)}$.
\dciteda{3.45 bzw. 3.47}
\end{proposition}

\begin{proof}
We will only prove the case $A \in \Adjlco{X}$.
Let $(T_t)_\indtGr$ be the $C_0$-group generated by $A$.
There exists $a_t \in \mIMls(X)$ such that $j(T_t x) = a_t \multis j(x)$ for all $x \in X$.
Define $S_t \colon \ternh{X} \to \ternh{X}, z \mapsto a_t \multis z$.
We conclude
\begin{equation}\label{normatz}
\norm{a_t \multis (j(x) \multis j(y)^* \multis j(z)) - j(x) \multis j(y)^* \multis j(z) }
\leq  \norm{ T_t x - x }\, \norm{j(y)^* \multis j(z)}  \to 0
\end{equation}
for $t \to 0$ and for all $x,y,z \in X$.
Because
\[ \cmlin\left\{ x_1 \multis x_2^* \multis x_3 \multis x_4^* \cdots x_{2n+1} \setfdg
                n \geq 0, x_1, \dots ,x_{2n+1} \in j(X) \right\} \]
is dense in $\ternh{X}$,
it follows with \cite[Proposition I.5.3]{EngelNagelSemigroups}
that $(S_t)_\indtGr$ is a \mbox{$C_0$-group}.

Let $B$ denote the generator of $(S_t)_\indtGr$.
We get
\[ j(Ax) = \lim_{t \to 0} \frac{a_t \multis j(x) - j(x)}{t} = B(j(x)) \]
for all $x \in X$.
The uniqueness of $B$ follows directly from the fact
that
\begin{align*}
W_{j(D(A)),j(X)} :=
   \cmlin &\left\{ x_1 \multis x_2^* \multis x_3 \multis x_4^* \cdots x_{2n+1}  \setfdg  \right. \\
          &\,\,\left. n \geq 1, x_1 \in j(D(A)), x_2 ,\dots, x_{2n+1} \in j(X) \right\}
\end{align*}
is a core for $B$, \cmie, $W_{j(D(A)),j(X)}$ is dense in $D(B)$ for the graph norm.
\end{proof}

Using the above proposition together with the Theorems \ref{zshgRegOpUnbeschrMult} and \ref{SatzLumerPhillipsOR},
we obtain the equivalence of (b) and (c) in Theorem \ref{thMainChar}.


Similar to regular operators, every unbounded multiplier has an adjoint:

\begin{theorem}
Let $A \in \Adjloru{X}$.
There exists a unique operator  $A^* \colon D(A^*) \subseteq X \to X$,
called the \defemph{adjoint} of $A$, such that
$\mri \cmpmatrix{ 0 & A \\ A^* & 0 } \in \Adjlco{C_2(X)}$.
We have: $A^* \in \Adjloru{X}$ and $(A^*)^* = A$.\dciteda{3.50, 3.60}
\end{theorem}

\begin{proof}
The existence follows directly from the definition of an
unbounded multiplier. To show the uniqueness, let $B_i \colon D(B_i)
\subseteq X \to X$ ($i \in \{1,2\}$) be such that $\tilde{B}_i
:= \mri \cmpmatrix{ 0 & A \\ B_i & 0} \in \Adjlco{C_2(X)}$.
Due to Proposition \ref{grHochliftenAufTX},
we find a regular operator 
$\tilde{C}_i \in \Multwor{\ternh{C_2(X)}}$ satisfying
$j \circ \tilde{B}_i = \tilde{C}_i \circ j\restring_{D(\tilde{B}_i)}$.
Let $(T_t^{(i)})_\indtGr$ be the $C_0$-group generated by $\tilde{B}_i$.
There exits $a_t^{(i)} \in \mIMls(X)$ such that
$j(T_t^{(i)} x) = a_t^{(i)} \multis j(x)$ for all $x \in X$.
Note that $\tilde{C}_i$ is the generator of the $C_0$-group
that is defined by 
$S_t^{(i)} \colon \ternh{X} \to \ternh{X}, z \mapsto a_t^{(i)} \multis z$.
For all $k \in \mbbN \cup \{0\}$, $x_1 \in D(\tilde{B}_i)$ and $x_2,\dots,x_{2k+1} \in C_2(X)$
we obtain with $z := j(x_2)^* \multis j(x_3) \multis j(x_4)^* \cdots j(x_{2k+1})$
\[ j(\tilde{B}_i(x_1)) \multis z
=  \lim_{t \to 0} \frac{a_t^{(i)} \multis j(x_1) - j(x_1)}{t}  \multis z
=  \tilde{C}_i \bigl( j(x_1) \multis z \bigr). \]

Set $W_i := W_{D(\tilde{B}_i),C_2(X)}$.
Using $\ternh{C_2(X)} \cong C_2(\ternh{X})$ (\cite[8.3.12.(4)]{BlecherLeMerdy04OpAlg}),
we conclude that $\tilde{C}_i\restring_{W_i}$ has the same form as $\tilde{B}_i$,
namely $\mri \cmpmatrix{ 0 & * \\ * & 0 }$.
Moreover, $W_i$ is a core for $\tilde{C}_i$.
Therefore, $\tilde{C}_i$ has the form $\mri \cmpmatrix{ 0 & * \\ * & 0 }$.
Using Theorem \ref{zshgRegOpUnbeschrMult},
there exists $\check{C}_i \in \Multwor{\ternh{X}}$
such that
\[ \tilde{C}_i = \mri \cmpmatrix{0 & \check{C}_i \\ \check{C}_i^* & 0}. \]
The set $W := W_{j(D(A)),j(X)}$
is a core for $\check{C}_1$ and for $\check{C}_2$.
Using $\check{C}_1\restring_W = \check{C}_2\restring_W$, we conclude
\[ \check{C}_1 = \overline{\check{C}_1\restring_W} = \overline{\check{C}_2\restring_W} = \check{C}_2. \]
It follows $\tilde{C}_1 = \tilde{C}_2$,
thus $\tilde{B}_1 = \tilde{B}_2$ and $B_1 = B_2$.

With $\check{C}_i^{**} = \check{C}_i$ (\cite[Corollary 9.4]{LanceHmod})
and $\check{C}_i^* \in \Multwor{\ternh{X}}$ (\cite[Corollary~9.6]{LanceHmod})
one can show that
\[ \hat{C}_i := \mri \cmpmatrix{0 & \check{C}_i^* \\ \check{C}_i & 0} \in \Multwor{\ternh{X} \oplus \ternh{X}} \]
is skew-adjoint.
Using $\ternh{X} \oplus \ternh{X} \cong C_2(\ternh{X})$,
we get $\hat{C}_i \in \Adjloru{C_2(\ternh{X})}$.
Set $Y := j_{C_2(X)}(C_2(X))$.
One can show that $\hat{C}_i\restring_Y \in \Adjloru{Y}$,
thus we obtain $\mri \cmpmatrix{0 & A^* \\ A & 0} \in \Adjloru{C_2(X)}$,
and the remaining assertions follow.
\end{proof}

\subsection{Connection to operators on Hilbert spaces}

In this section we characterize unbounded skew-adjoint multipliers
using $C_0$-groups on Hilbert spaces.

\begin{theorem}\label{charAdjlcoHR}
Let $(\pi,K,K_1,K_2)$ be a $\ternh{X}$-embedding of $X$
and $\eta := \pi \circ j \colon X \to \mLinStet(K)$.
Let $A \colon D(A) \subseteq X \to X$ be linear.
The following assertions are equivalent:
\begin{enumaequiv}
\item $A \in \Adjlco{X}$.

\item $A$ is densely defined with $\rho(A) \neq \emptyset$, and
there exists a $C_0$-group $(b_t)_\indtHG$
of unitary elements in $\mLinStet(K_2)$ such that:
\begin{enumaufz}
\item $L_{b_t} \colon \eta(X) \to \eta(X), y \mapsto b_t y$
and $L_{b^*_t}$ define $C_0$-groups on $\eta(X)$ and
\item for the generator $C$ of $(L_{b_t})_\indtGr$, we have $j(D(A)) \subseteq D(C)$.


\dciteda{3.66}%
\end{enumaufz}
\end{enumaequiv}
\end{theorem}

\begin{proof}
At first we show that (a) implies (b).
Let $(T_t)_\indtGr$ denote the $C_0$-group generated by $A$.
There exists $a_t \in \mIMls(X)$ such that $j(T_t x) = a_t \multis j(x)$.
Define $b_t := \pi(a_t)\restring_{K_2}$.
Because $\erzl{\pi(\ternh{X})K_1}$ is dense in $K_2$,
we have $b_t(K_2) \subseteq K_2$.
With equation \eqref{normatz} it follows for all $z \in \ternh{X}$ and $\xi \in K_1$
\[
   \norm{ b_t(\pi(z)(\xi)) - \pi(z)(\xi) }
=  \norm{ (\pi(a_t) \pi(z) - \pi(z))(\xi) }
\leq  \norm{ a_t \multis z - z } \, \norm{\xi}
\to  0 \]
for $t \to 0$.
Therefore $(b_t)_\indtGr$ is a $C_0$-group on $K_2$.
Using the theorem of Hille-Yosida, we obtain
that $A$ is densely defined with $\rho(A) \neq \emptyset$.
It is straightforward to show (i) and (ii).

To prove the other direction, define $d_t \colon K \to K, \xi \mapsto b_t(p_{K_2}(\xi))$
where $p_{K_2}$ denotes the projection from $K$ onto $K_2$
and $S_t := L_{d_t} \colon \eta(X) \to \eta(X)$.
Then it can be checked that $T_t := \eta^{-1} \circ S_t \circ \eta$
defines a $C_0$-group on $X$ with generator $A$.
\dremark{Bew. ist sehr knapp.}%
\end{proof}

With the above theorem, we obtain the equivalence of (a) and (d) in Theorem \ref{thMainChar}.
So we finished the proof of this theorem.



\begin{proposition}\label{liftAdjlcoToHS}
Let $A \in \Adjlco{X}$.
There exist a Hilbert space $K$, a completely isometric map $\eta \colon X \to \mLinStet(K)$
and a generator $C$ of a unitary $C_0$-group on $K$ such that
$\eta(Ax) = C \circ \eta(x)$ for all $x \in D(A)$.
\dremark{Neu!}
\end{proposition}

\begin{proof}
We use the notation of Theorem \ref{charAdjlcoHR}.
Using this result, we find a $C_0$-group $(b_t)_\indtGr$
with the properties stated in Theorem \ref{charAdjlcoHR}.
Let $\tilde{C}$ be the generator of this $C_0$-group.
With Theorem \ref{charAdjlcoHR} we get
\[ \eta(Ax)(\xi)
=  \lim_{t \to 0} \frac{(b_t \circ \eta(x))(\xi) - \eta(x)(\xi)}{t}
=  \tilde{C}(\eta(x)(\xi))  \]
for all $x \in D(A)$ and $\xi \in K$,
thus $\eta(Ax) = \tilde{C} \circ \eta(x)$.
The assertion follows from setting $C := \tilde{C} \oplus 0_{K_2^\perp}$.
\end{proof}

\subsection{The strict $X$-topology}

In order to transform a $C_0$-group $(b_t)_\indtGr$ from a Hilbert space to an operator space,
we have to impose a stronger condition than strong continuity
on $(b_t)_\indtGr$ (cf. Theorem \ref{charAdjlcoHR}), namely
\[ \lim_{t \to 0} b_t \circ y = y \quad\text{for all } y \in \eta(X). \]
In order to formulate this in a more abstract fashion,
we introduce the following topology for an operator space $X \subseteq \mLinStet(H_0,H)$:

\begin{definitn}
For all $x \in X$ a seminorm is defined by
\[ p_x \colon L(H) \to \mbbR, T \mapsto \norm{T \circ x}. \]
The locally convex topology on $\mLinStet(H)$, generated by $(p_x)_{x \in X}$,
is called the \defemph{strict $X$-topology}.
\end{definitn}

The norm topology on $\mLinStet(H)$ is finer than
the strict $X$-topology on $\mLinStet(H)$.
If $X$ is unital and $H_0 = H$, these two topologies coincide.

The proof of the following proposition is straightforward:

\begin{proposition}
Let $\erzl{X H_0}$ be dense in $H$.
The strict $X$-topology on $\mLinStet(H)$ is finer than the strong topology.
\dciteda{3.70}%
\end{proposition}

It is not hard to see
that the assumption in the above proposition can always be fulfilled
by choosing an appropriate embedding of $X$.
\dremark{Sei $X \subseteq \mLinStet(H_0,K)$.
Setze $H := \overline{\erzl{X H_0}}$.}
\pagebreak 

\begin{proposition}
Let $K$ be a Hilbert space and
$\eta \colon X \to \mLinStet(K)$ a completely isometric map.
Let $C \colon D(C) \subseteq K \to K$ be the generator of a
$C_0$-group $(b_t)_\indtGr$ of unitary elements in $\mLinStet(K)$ such that
\begin{enumaufz}
\item $(b_t)_t$ converges in the strict $\eta(X)$-topology to $\Id_K$ for $t \to 0$ and
\item $b_t \eta(X), b_t^* \eta(X) \subseteq \eta(X)$
for all $\indtGr$.
\end{enumaufz}
Then there exists $A \in \Adjlco{X}$ such that
\dremark{$\eta(D(A))(K) \subseteq D(C)$ and}%
\dremark{Neu!}%
\[ \eta(Ax)(\xi) = C(\eta(x)(\xi)) \quad\text{for all } x \in D(A), \xi \in K. \]
\end{proposition}

\begin{proof}
The proof is similar to the proof of Proposition \ref{liftAdjlcoToHS}.
\end{proof}


\subsection{Pertubation theory}

The following result generalizes a pertubation result of Damaville
(\cite[Proposition 2.1]{Damaville04RegulariteDesOp}, see also \cite{Damaville07RegulariteDOpReg})
from Hilbert \cstern{}modules to operator spaces:

\begin{theorem}
\begin{enumaufz}
\item If $A \in \Multlco{X}$ and $B \in \Multlor{X}$ then $A + B \in \Multlco{X}$.

\item If $A \in \Adjlco{X}$ and $B \in \Adjlor{X}$ with $B^* = -B$
then $A+B \in \Adjlco{X}$.

\item If $A \in \Adjloru{X}$ and $B \in \Adjlor{X}$
then $A+B \in \Adjloru{X}$.
\end{enumaufz}
\end{theorem}

\begin{proof}
Part (i) follows directly with \cite[Theorem III.1.10]{EngelNagelSemigroups}.

To prove (ii), we notice that there exists
$\check{A} \in \Adjlco{\ternh{X}} = \Multwor{\ternh{X}}$
satisfying $j \circ A = \check{A} \circ j\restring_{D(A)}$
due to Theorem \ref{grHochliftenAufTX}.
There exists $\check{B} \in \Adjlor{\ternh{X}} \cong \Adjhm{\ternh{X}}$ skew-adjoint
satisfying $j \circ B = \check{B} \circ j$.
Using \cite[Proposition 2.1]{Damaville04RegulariteDesOp}, we conclude
$\check{A} + \check{B} \in \Adjlco{\ternh{X}}$.
One can show that $A + B \in \Adjlco{X}$.

Part (iii) is a consequence of part (ii).
\end{proof}



\bibliographystyle{amsalpha}  
\bibliography{Journal}






\end{document}